\documentclass{article}
\usepackage{amssymb}
\usepackage{amsmath,amsthm}
\usepackage{amsfonts}
\usepackage{graphicx}
\usepackage{graphics}

\vspace{.9cm}
\newtheorem{Theorem}{Theorem}[section]

\newtheorem{cor}{Corollary}[section]
\newtheorem{lemma}{Lemma}[section]

\begin{document}

\date{\small\textsl{\today}}
\title{ Convergence of general composite iterative method for infinite
family of nonexpansive mappings in Hilbert spaces}
\author{  \large Vahid Darvish
          and ~\large S. M. Vaezpour
\footnote{Corresponding author.\newline {\em  E-mail addresses:}
\newline vahiddarvish@aut.ac.ir (V. Darvish),
\newline vaez@aut.ac.ir (S. M. Vaezpour).}
           $\vspace{.2cm} $  \\
\small{\em Department of Mathematics and Computer Science,}\vspace{-1mm}\\
\small{\em Amirkabir University of Technology,}\vspace{-1mm}\\
\small{\em Hafez Ave., P.O. Box 15875-4413, Tehran, Iran }\\}
\maketitle

\begin{abstract} In this paper by using $W_{n}$-mapping, we introduce
a composite iterative method for finding a common fixed point for
infinite family of nonexpansive mappings and a solution of a certain
variational inequality. Furthermore, the strong convergence of the
proposed iterative method is established. Finally, some simulation
examples are presented. Our results improve and extend the previous
results.

\medskip


\medskip


\end{abstract}

\section{Introduction and preliminaries}
Let $H$ be a real Hilbert space, whose inner product and norm are
denoted by $\langle.,.\rangle$ and $\|. \|$, respectively. Let $C$
be a nonempty closed convex subset of $H$ and $T$ is a nonlinear
mapping. We use $F(T)$ to denote the set of fixed points of $T$
(i.e., $F(T)=\{x\in H:Tx=x\}$). Recall that a self mapping $T$ of
$C$ is {\it nonexpansive } if $\|Tx-Ty\|\leq\|x-y\|, \forall x,y\in
H$ and is a contraction, if there exists a constant $\alpha \in
(0,1)$ such that $\|
Tx-Ty\|\leq\alpha\|x-y\|$ for all $x,y\in C$.\\
A bounded linear operator $A$ on $H$ is called {\it strongly
positive } with coefficient $ \bar{\gamma}>0 $ if,
$$ \langle Ax,x \rangle \geq \bar{\gamma} \|x\|^2, \forall x\in H  .$$

\par
\

In 2005,  Kim and Xu \cite{kim} introduced the following iteration process:\\

\begin{eqnarray}
x_{0} &=& x\in C  \  chosen \ arbitrary \ \nonumber  ,\\
y_{n} &=&  \beta_{n}x_{n}+(1-\beta_{n})Tx_{n},\nonumber \\
x_{n+1} &=&  \alpha_{n}u+(1-\alpha_{n})y_{n} \nonumber \label{o}.\\
\end{eqnarray}

\noindent They proved in a uniformly smooth Banach space, the
sequence $\{x_{n}\}$ defined by (\ref{o}) converges strongly to a
fixed point of $T$.
In 2009 Cho and Qin \cite{yeol} considered the following composite iterative algorithm :\\

$$x_{0}\in H  \ chosen \ arbitrary  ,$$
$$z_{n}=\gamma_{n}x_{n}+(1-\gamma_{n})Tx_{n},$$
$$y_{n}=\beta_{n}x_{n}+(1-\beta_{n})Tz_{n}, $$
\begin{equation}\label{p}
\ x_{n+1}=\alpha_{n}\gamma f(x_{n})+\delta_{n}x_{n}+((1-\delta_{n})I-\alpha_{n}A)y_{n}, \ \forall \ n\geq0 \ . \nonumber
\end{equation}
\

In 2009 Wangkeeree and Kamraksa \cite{wan} introduced a new iterative scheme:\\
$$ \ x_{0} =x\in C  \ chosen  \ arbitrary ,$$
$$z_{n}=\gamma_{n}x_{n}+(1-\gamma_{n})W_{n}x_{n},$$
$$y_{n}=\beta_{n}x_{n}+(1-\beta_{n})W_{n}z_{n},$$
\begin{equation}\label{wan}
\ x_{n+1}=\alpha_{n}\gamma
f(x_{n})+\delta_{n}x_{n}+((1-\delta_{n})I-\alpha_{n}A)P_{C}(y_{n}-\lambda_{n}By_{n}),
\end{equation}
\par

\noindent where the mapping $W_{n}$ defined by Shimoji and Takahashi
\cite{shi}, as follows:
\begin{eqnarray}
U_{n,n+1}&=&I,\nonumber\\
U_{n,n}&=&\gamma_{n}T_{n}U_{n,n+1}+(1-\gamma_{n})I,\nonumber\\
U_{n,n-1}&=&\gamma_{n-1}T_{n-1}U_{n,n}+(1-\gamma_{n-1})I,\nonumber\\
&\vdots&\nonumber\\
U_{n,k}&=&\gamma_{k}T_{k}U_{n,k+1}+(1-\gamma_{k})I,\nonumber\\
U_{n,k-1}&=&\gamma_{k-1}T_{k-1}U_{n,k}+(1-\gamma_{k-1})I,\nonumber\\
&\vdots&\nonumber\\
U_{n,2}&=&\gamma_{2}T_{2}U_{n,3}+(1-\gamma_{2})I,\nonumber\\
W_{n}&=&U_{n,1}=\gamma_{1}T_{1}U_{n,2}+(1-\gamma_{1})I,\nonumber\label{12}\\
\end{eqnarray}
where $\gamma_{1},\gamma_{2},...$ are real numbers such that $0\leq
\gamma_{n}\leq1, T_{1},T_{2},...$ are an infinite family of mappings
of $H$ into itself, note that the nonexpansivity of each $T_{i}$
ensures the nonexpansivity of $W_{n}$. \\

In 2010 Singthong and Suantai \cite{sing} introduced an iterative
method as follows:
\begin{eqnarray}
x_{0}&=&x\in C \ chosen \ arbitrary ,\nonumber\\
y_{n}&=&\beta_{n}x_{n}+(1-\beta_{n})K_{n}x_{n} ,\nonumber\\
x_{n+1}&=&P_{C}(\alpha_{n}\gamma f(x_{n})+(I-\alpha_{n}A)y_{n})
,\nonumber\label{va}\\
\end{eqnarray}
where $K$-mapping defined by Kangtunyakarn and Suantai \cite{suan},
as follows:\\

\begin{eqnarray}
U_{n,1}&=&\lambda_{n,1}T_{1}+(1-\lambda_{n,1})I,\nonumber\\
U_{n,2}&=&\lambda_{n,2}T_{2}U_{n,1}+(1-\lambda_{n,2})U_{n,1},\nonumber\\
U_{n,3}&=&\lambda_{n,3}T_{3}U_{n,2}+(1-\lambda_{n,3})U_{n,2},\nonumber\\
&\vdots&\nonumber\\
U_{n,N-1}&=&\lambda_{n,N-1}T_{N-1}U_{n,N-1}+(1-\lambda_{n,N})U_{n,N-1},\nonumber\\
K_{n}=U_{n,N}&=&\lambda_{n,N}T_{N}U_{n,N-1}+(1-\lambda_{n,N})U_{n,N-1} \ ,\nonumber\label{64}\\
\end{eqnarray}
where $\{T_{i}\}_{i=1}^N$ are finite family of nonexpansive mappings
and the sequences $\{\lambda_{n,i}\}_{i}^N$ are in $[0,1]$. The mapping $K_{n}$ is called the $K$-mapping generated by  $T_{1},T_{2},\ldots,T_{N}$ and $\lambda_{n,1},\lambda_{n,2},\ldots,\lambda_{n,N}$.\\

Through out this paper inspired by Singthong and Suantai \cite{sing}
and Wangkeeree and Kamraksa \cite{wan},
we introduce a composite iteration method for infinite family of nonexpansive mappings as follows :\\
\par

$$x_{0}=x\in C \ chosen \ arbitrary  ,$$
$$z_{n}=\gamma_{n}x_{n}+(1-\gamma_{n})W_{n}x_{n},$$
$$y_{n}=\beta_{n}x_{n}+(1-\beta_{n})W_{n}z_{n},$$
\begin{equation}\label{a}
\ x_{n+1}=P_{C}[\alpha_{n}\gamma
f(x_{n})+\delta_{n}x_{n}+((1-\delta_{n})I-\alpha_{n}A)y_{n}],
\end{equation}
\par
\noindent where $W_{n}$ is defined by (\ref{12}), $f$ is a
contraction on $H$, $A$ is a strongly positive linear bounded
self-adjoint operator with
the coefficient $\bar{\gamma}>0$ and $0< \gamma < \frac{\bar{\gamma}}{\alpha}$.\\
Then by using this iteration we prove the existence of a common
fixed point for infinite family of nonexpansive mappings and the
solution of a certain variational inequality.

\

We need the following lemmas for the proof of our main results.\\

\begin{lemma}\label{r} The following inequality holds in a Hilbert space
H,
$$\|x+y\|^2 \leq \|x\|^2 + 2\langle y,x+y\rangle, \forall x,y \in H \ .$$
\end{lemma}
\begin{lemma}\label{1} \cite{xu} Assume $\{\alpha_n\}$ is a sequence of nonnegative real
numbers such that $\alpha_{n+1} \leq (1-\gamma_n)\alpha_n + \delta_n
\ \ n\geq 1$ , where $\{\alpha_n \}$ is a sequence in $(0,1)$ and
$\delta_n$ is a
sequence in $\mathbb{R}$ such that :\\
\begin{enumerate}
\item $\sum_{n=1}^{\infty}\gamma_{n}=\infty ,$
\item $\limsup_{n \to \infty}(\frac{\delta_{n}}{\gamma_{n}})\leq 0
$ or $\sum_{n=1}^{\infty}|\delta_{n}|< \infty , $
\end{enumerate}
then $\lim_{n\to\infty}\alpha_{n}=0 .$
\end{lemma}
\begin{lemma}\label{d} \cite{marino} Assume that A is a strongly positive linear bounded
self-adjoint operator on a Hilbert space $H$ with coefficient
$\bar{\gamma}$ and $0<\rho\leq \|A\|^{-1}$, then $\|I-\rho
A\|\leq1-\rho\bar{\gamma}$.
\end{lemma}
\begin{lemma}\cite{shi} Let $C$ be nonempty closed convex subset of a Hilbert
space, let $T_{i}:C\longrightarrow C$ be an infinite family of
nonexpansive mappings with $\bigcap_{i=1}^{\infty}F(T_{i})\neq
\emptyset$ and let ${\gamma_{i}}$ be a real sequence such that
$0<\gamma_{i}\leq\gamma<1$ for all $i\geq1$ then ,
\begin{enumerate}
\item $W_{n}$ is nonexpansive and $F(W_{n})=\bigcap_{i=1}^{n}F(T_{i})$ for
each $n\geq 1$ ,
\item For each $x\in C$ and for each positive integer $k$, the
$\lim_{n\to\infty}U_{n,k}$ exists ,
\item The mapping $W:C\longrightarrow C$ defined by :
$$Wx:=\lim_{n\to\infty}W_{n}x=\lim_{n\to\infty}U_{n,1}x \  \ \ x\in
C,$$ \\

is a nonexpansive mapping satisfying
$F(W)=\bigcap_{i=1}^{\infty}F(T_{i})$ and is called the $W$-mapping
generated by $T_{1},T_{2},...$ and $\gamma_{1},\gamma_{2},... \ .$
\end{enumerate}
\end{lemma}
\begin{lemma}\cite{shi} Let $C$ be a nonempty closed convex subset of a Hilbert
space $H$, let ${T_{i}:C\longrightarrow C }$ be an infinite family
of nonexpansive mappings with $\bigcap_{i=1}^{\infty}F(T_{i})\neq
\emptyset $  and let ${\gamma_{i}}$ be a real sequence such that
$0<\gamma_{i}\leq \gamma<1$ for all $i\geq1$, \ if $K$ is any
bounded subset of $C$ then,
$$\limsup_{n \to\infty}\|Wx-W_{n}x\|=0 \ \ x\in K.$$
\end{lemma}
\begin{lemma}\label{e}\cite{marino} Let $H$ be a Hilbert space, let $A$ be a strongly positive
linear bounded self-adjoint operator with coefficient
$\bar{\gamma}>0$. Assume that $ 0 < \gamma <
\frac{\bar{\gamma}}{\alpha}$, let $T$ be a nonexpansive mapping with
a fixed point $x_{t}$ of the contraction, $$x\longmapsto t\gamma
f(x)+(I-tA)Tx,$$ then ${x_{t}}$ converges strongly as $ t
\rightarrow 0$ to a fixed point $\bar{x}$ of $T$ which solves the
variational inequality $\langle(A-\gamma f)\bar{x},\bar{x}-z\rangle
\leq0 \ \forall z\in F(T).$
\end{lemma}
\begin{lemma}\cite{suan}
Let $C$ be a nonempty closed convex subset of strictly convex Banach
space. Let $\{T_{i}\}_{i=1}^N$ be a finite family of nonexpansive
mappings of $C$ into itself with $\bigcap_{i=1}^N
F(T_{i})\neq\emptyset$, and let $\lambda_{1},\ldots,\lambda_{N}$ be
real numbers such that $0<\lambda_{i}<1$ for every $i=1,\ldots,N-1$
and $0<\lambda_{N}\leq1$. Let $K$ be the $K$-mapping of $C$ into itself generated by $T_{1},\ldots,T_{N}$ and $\lambda_{1},\ldots,\lambda_{N}$. Then, \\
\begin{equation}\label{5}
\ \ \ \ \  F(K)=\bigcap_{i=1}^N F(T_{i}).
\end{equation}
\end{lemma}
\begin{lemma}\cite{sing}
Let $C$ be a nonempty closed convex subset of a Banach space. Let
$\{T_{i}\}_{i=1}^N$ be a finite family of nonexpansive mappings of
$C$ into itself and $\{\lambda_{n,i}\}_{i=1}^N$ sequences in $[0,1]$
such that $\lambda_{n,i}\to\lambda_{i}$, as $n\to\infty$ ,
$(i=1,2,\ldots,N)$. Moreover, for every $n\in \mathbb{N}$ , K and
$K_{n}$ be the $K-mapping$ generated by $T_{1},\ldots,T_{N}$ and
$\lambda_{1},\ldots,\lambda_{N}$ and $T_{1},\ldots,T_{N}$ and
$\lambda_{n,1},\lambda_{n,2},\ldots,\lambda_{n,N}$,
respectively.\\
Then, for every bounded sequence $x_{n}\in C$, we have
$\lim_{n\to\infty}\|K_{n}x_{n}-Kx_{n}\|=0$.
\end{lemma}

\section{Main Results}
In this section, we prove strong convergence of the sequences
$\{x_{n}\}$ defined by the iteration scheme (\ref{a}), for finding a
common fixed point of infinite family of nonexpansive mappings which
solves the variational inequality.\\

\begin{Theorem} Let $C$ be a closed convex subset of a real Hilbert
space $H$. Let $f$ be a contraction of $C$ into itself, let $A$ be a
strongly positive linear bounded operator with coefficient
$\bar{\gamma}>0$ and $\{T_{i}:C \longrightarrow C\}$ be an infinite
family of nonexpansive mappings. Assume that
$0<\gamma<\frac{\bar{\gamma}}{\alpha}$ and
$F=\bigcap_{i=1}^{\infty}F(T_{i})\neq\emptyset$. Let $x_{0} \in C$,
given that $\{\alpha_{n}\}, \{\beta_{n}\},\{\gamma_{n}\}\ and \
\{\delta_{n}\}$ be sequences in $[0,1]$ satisfying the following
conditions :\\
$C_{1}: \lim_{n \rightarrow \infty}\alpha_{n}=0 \,
\sum_{n=1}^\infty\alpha_{n}=\infty$,\\
$C_{2} : 0<\liminf_{n
\rightarrow\infty}\delta_{n}\leq\limsup_{n\rightarrow\infty}\delta_{n}<1$,\\
$C_{3}\colon \sum_{n=1}^{\infty}|\gamma_{n}-\gamma_{n-1}|<\infty$,\\
$C_{4}\colon \sum_{n=1}^{\infty}|\alpha_{n}-\alpha_{n-1}|<\infty$,\\
$C_{5}\colon \sum_{n=1}^{\infty}|\beta_{n}-\beta_{n-1}|<\infty$,\\
$C_{6}\colon (1+\beta_{n})\gamma_{n}-2\beta_{n}>d$ for some $\quad  d\in (0,1)$,\\
then the sequence $\{x_{n}\}$ defined by (\ref{a}) converges
strongly to $q\in F$ which solves the variational inequality
$\langle\gamma
f(q)-Aq,p-q\rangle\leq0$, $\forall p\in F$.\\
\end{Theorem}
\noindent{proof}\ : Since $\alpha_{n} \to 0$ as $n\to \infty$
without loss of generality we have \\
$\alpha_{n}<(1-\delta_{n})\|A\|^{-1}\ \ \forall n\geq0$, noticing
that $A$ is a bounded linear self-adjoint operator with,
$$\|A\|=\sup \{|<Ax,x>|\ :x\in H , \|x\|=1\},$$
we have,\\
\begin{eqnarray*}
<((1-\delta_{n})I-\alpha_{n}A)x,x>&=&(1-\delta_{n})<x,x>-\alpha_{n}<Ax,x>  \\
& \geq & (1-\delta_{n})-\alpha_{n}\|A\|\geq0  , \\
\end{eqnarray*}
then $(1-\delta_{n})I-\alpha_{n}A$ is positive. Also,\\
\begin{eqnarray}
\|(1-\delta_{n})I-\alpha_{n}A\|&=&\sup\{|<((1-\delta_{n})I-\alpha_{n}A))x,x>|,x\in
H,\|x\|=1\} \nonumber\\
&=&\sup\{1-\delta_{n}-\alpha_{n}<Ax,x>,x\in H , \|x\|=1\}\ \nonumber\\
&\leq& 1-\delta_{n}-\alpha_{n}\bar{\gamma} \nonumber \label{b}  .\\
\end{eqnarray}

\noindent Next we prove that $\{x_{n}\}$ is bounded. We pick $p\in F=\bigcap_{i=1}^{\infty}F(T_{i})=F(W)=F(W_{n}) .$\\
\begin{eqnarray*}
\|z_{n}-p\|&=&\|\gamma_{n}x_{n}+(1-\gamma_{n})W_{n}x_{n}-p\|\\
&=&\|\gamma_{n}(x_{n}-p)+(1-\gamma_{n})(W_{n}x_{n}-W_{n}p)\| \\
&\leq& \gamma_{n}\|(x_{n}-p)\|+(1-\gamma_{n})\|(x_{n}-p)\|\\
&=&\|x_{n}-p\| ,\\
\end{eqnarray*}
\noindent and we have,\\
\begin{eqnarray*}
\|y_{n}-p\|&=&\|\beta_{n}x_{n}+(1-\beta_{n})W_{n}z_{n}-p\|\\
&=&\|\beta_{n}(x_{n}-p)+(1-\beta_{n})(W_{n}z_{n}-W_{n}p)\|\\
&\leq& \beta_{n}\|x_{n}-p\|+(1-\beta_{n})\|z_{n}-p\|\\
&\leq&\beta_{n}\|x_{n}-p\|+(1-\beta_{n})\|x_{n}-p\|\\
&=&\|x_{n}-p\| .\\
\end{eqnarray*}
It follows that,\\
\begin{eqnarray*}
\|x_{n+1}-p\|&=&\|P_{C}[\alpha_{n}\gamma f(x_{n})+\delta_{n}x_{n}+((1-\delta_{n})I-\alpha
A)y_{n}]-P_{C}(p)\|\\
&\leq&\|\alpha_{n}\gamma f(x_{n})+\delta_{n}x_{n}+((1-\delta_{n})I-\alpha_{n}A)y_{n}-p\|\\
&=&\|\alpha_{n}(\gamma f(x_{n}-Ap)+\delta_{n}(x_{n}-p)+((1-\delta_{n})I-\alpha_{n}A)(y_{n}-p)\| ,\\
\end{eqnarray*}

\noindent by (\ref{b}) we have,\\
\begin{eqnarray*}
&\leq &\alpha_{n}\|\gamma f(x_{n})-Ap\|+\delta_{n}\|x_{n}-p\|+(1-\delta_{n}-\alpha_{n}\bar{\gamma})\|y_{n}-p\|\\
&\leq&\alpha_{n}\gamma\|f(x_{n})-f(p)\|+\alpha_{n}\|\gamma f(p)-Ap\|+\delta_{n}\|x_{n}-p\|+(1-\delta_{n}-\alpha_{n}\bar{\gamma})\|x_{n}-p\|\\
&\leq& \alpha_{n}\gamma\alpha\|x_{n}-p\|+\alpha_{n}\|\gamma f(p)-Ap\|+(1-\alpha_{n}\bar{\gamma})\|x_{n}-p\|\\
&=&[1-\alpha_{n}(\bar{\gamma}-\gamma\alpha)]\|x_{n}-p\|+\alpha_{n}\|\gamma f(p)-Ap\| .\\
\end{eqnarray*}

\noindent By simple induction we have $\|x_{n}-p\|\leq
max\{\|x_{0}-p\|,\frac{\|Ap-\gamma
f(p)\|}{\bar{\gamma}-\gamma\alpha}\}$, which  gives that the
sequence $\{x_{n}\}$ is bounded so are $\{y_{n}\}$ and
$\{z_{n}\}$.\\

\noindent Next we claim that,
$\lim_{n\to\infty}\|x_{n+1}-x_{n}\|=0$ .\\
We know that,\\
\begin{eqnarray*}
z_{n}&=&\gamma_{n}x_{n}+(1-\gamma_{n})W_{n}x_{n}\ ,\\
z_{n-1}&=&\gamma_{n-1}x_{n-1}+(1-\gamma_{n-1})W_{n-1}x_{n-1}.\\
\end{eqnarray*}
So we obtain,\\
$$z_{n}-z_{n-1}=(1-\gamma_{n})(W_{n}x_{n}-W_{n-1}x_{n-1})+\gamma_{n}(x_{n}-x_{n-1})+(\gamma_{n}-\gamma_{n-1})(x_{n-1}-W_{n-1}x_{n-1}) .$$\\
This implies that,\\
\begin{eqnarray*}
\|z_{n}-z_{n-1}\|&\leq&(1-\gamma_{n})\|W_{n}x_{n}-W_{n-1}x_{n-1}\|+\gamma_{n}\|x_{n}-x_{n-1}\|+|\gamma_{n}-\gamma_{n-1}|\|x_{n-1}-W_{n-1}x_{n-1}\|\\
                 &=&(1-\gamma_{n})\|W_{n}x_{n}-W_{n}x_{n-1}+W_{n}x_{n-1}-W_{n-1}x_{n-1}\|\\
                 &&+\gamma_{n}\|x_{n}-x_{n-1}\|+|\gamma_{n}-\gamma_{n-1}|\|x_{n-1}-W_{n-1}x_{n-1}\|\\
                 &\leq&(1-\gamma_{n})\|W_{n}x_{n}-W_{n}x_{n-1}\|+(1-\gamma_{n})\|W_{n}x_{n-1}-W_{n-1}x_{n-1}\|\\
                &&+\gamma_{n}\|x_{n}-x_{n-1}\|+|\gamma_{n}-\gamma_{n-1}|\|x_{n-1}-W_{n-1}x_{n-1}\| .\\
\end{eqnarray*}

\noindent On the other hand we have:\\
\begin{eqnarray*}
\|W_{n}x_{n-1}-W_{n-1}x_{n-1}\|&=&\|\gamma_{1}T_{1}U_{n,2}x_{n-1}-\gamma_{1}T_{1}U_{n-1,2}x_{n-1}\|\\
&\leq&\gamma_{1}\|U_{n,2}x_{n-1}-U_{n-1,2}x_{n-1}\|\\&=&\gamma_{1}\|\gamma_{2}T_{2}U_{n,3}x_{n-1}-\gamma_{2}T_{2}U_{n-1,3}x_{n-1}\|\\
&\leq&\gamma_{1}\gamma_{2}\|U_{n,3}x_{n-1}-U_{n-1,3}x_{n-1}\|\\
&\vdots&\\
&\leq&\gamma_{1}\gamma_{2}\ldots\gamma_{n-1}\|U_{n,n}x_{n-1}-U_{n-1,n}x_{n-1}\|\\
&\leq&M_{1}\prod_{i=1}^{n-1}\gamma_{i},\\
\end{eqnarray*}
\begin{equation}\label{36} \end{equation}
where $M_{1}\geq0$ is an appropriate constant such that:\\
$$\|U_{n,n}x_{n-1}-U_{n-1,n}x_{n-1}\|\leq M_{1} \ \ \forall  \
n\geq0.$$ \noindent Note that the boundedness of $x_{n}$ and the
nonexpansivity of $T_{n}$ ensure the existence of $M_{1}$.
So we have,\\
\begin{eqnarray*}
\|z_{n}-z_{n-1}\|&\leq&\gamma_{n}\|x_{n}-x_{n-1}\|+(1-\gamma_{n})M_{1}\prod_{i=1}^{n-1}\gamma_{i}\\
&&+(1-\gamma_{n})\|x_{n}-x_{n-1}\|+|\gamma_{n}-\gamma_{n-1}|\|x_{n-1}-W_{n-1}x_{n-1}\|\\
&=&\|x_{n}-x_{n-1}\|+(1-\gamma_{n})M_{1}\prod_{i=1}^{n-1}\gamma_{i}+|\gamma_{n}-\gamma_{n-1}|\|x_{n-1}-W_{n-1}x_{n-1}\| .\\
\end{eqnarray*}

\noindent Similar to (\ref{36}), we have
$$\|U_{n,n}z_{n-1}-U_{n-1,n}z_{n-1}\|\leq M_{2}.$$
So,
\begin{eqnarray*}
\|y_{n}-y_{n-1}\|&=&\|\beta_{n}x_{n}+(1-\beta_{n})W_{n}z_{n}-\beta_{n-1}x_{n-1}-(1-\beta_{n-1})W_{n-1}z_{n-1}\|\\
&=&\|\beta_{n}x_{n}-\beta_{n}x_{n-1}+\beta_{n}x_{n-1}-\beta_{n-1}x_{n-1}+(1-\beta_{n})(W_{n}z_{n}-W_{n}z_{n-1})\\
&&+(1-\beta_{n})(W_{n}z_{n-1}-W_{n-1}z_{n-1})+(1-\beta_{n})W_{n-1}z_{n-1}-(1-\beta_{n-1})W_{n-1}z_{n-1}\|\\
&\leq&\|\beta_{n}\left(x_{n}-x_{n-1}\right)+(\beta_{n}-\beta_{n-1})x_{n-1}+(1-\beta_{n})(W_{n}z_{n}-W_{n}z_{n-1})\\
&&+(1-\beta_{n})(W_{n}z_{n-1}-W_{n-1}z_{n-1})+(1-\beta_{n})W_{n-1}z_{n-1}-(1-\beta_{n-1})W_{n-1}z_{n-1}\|\\
&\leq&\beta_{n}\|x_{n}-x_{n-1}\|+(1-\beta_{n})\|z_{n}-z_{n-1}\|+(1-\beta_{n})M_{2}\prod_{i=1}^{n-1}\gamma_{i}\\
&&+|\beta_{n}-\beta_{n-1}|\|x_{n-1}-W_{n-1}z_{n-1}\|\\
&\leq& \beta_{n}\|x_{n}-x_{n-1}\|+(1-\beta_{n})\|x_{n}-x_{n-1}\|+(1-\beta_{n})|\gamma_{n}-\gamma_{n-1}|\|x_{n-1}-W_{n-1}x_{n-1}\|\\
&&+(1-\beta_{n})(1-\gamma_{n})M_{1}\prod_{i=1}^{n-1}\gamma_{i}+(1-\beta_{n})M_{2}\prod_{i=1}^{n-1}\gamma_{i}\\
&&+|\beta_{n}-\beta_{n-1}|\|x_{n-1}-W_{n-1}z_{n-1}\|\\
&=&\|x_{n}-x_{n-1}\|+(1-\beta_{n})|\gamma_{n}-\gamma_{n-1}|\|x_{n-1}-W_{n-1}x_{n-1}\|\\
&&+(1-\beta_{n})(1-\gamma_{n})M_{1}\prod_{i=1}^{n-1}\gamma_{i}+(1-\beta_{n})M_{2}\prod_{i=1}^{n-1}\gamma_{i}\\
&&+|\beta_{n}-\beta_{n-1}|\|x_{n-1}-W_{n-1}z_{n-1}\| .\\
\end{eqnarray*}
Therefore,
\begin{eqnarray*}
 \|x_{n+1}-x_{n}\|&=&\|P_{C}[\alpha_{n}\gamma
f(x_{n})+\delta_{n}x_{n}+((1-\delta_{n})I-\alpha_{n}A)y_{n}]
-P_{C}[\alpha_{n-1}\gamma f(x_{n-1})\\
&&+\delta_{n-1}x_{n-1}+((1-\delta_{n-1})I-\alpha_{n-1}A)y_{n-1}]\|\\
&\leq&\|\alpha_{n}\gamma
f(x_{n})+\delta_{n}x_{n}+((1-\delta_{n})I-\alpha_{n}A)y_{n}
-\alpha_{n-1}\gamma f(x_{n-1})-\delta_{n-1}x_{n-1}\\
&&-((1-\delta_{n-1})I-\alpha_{n-1}A)y_{n-1}\|\\
&\leq&\|((1-\delta_{n})I-\alpha_{n}A)(y_{n}-y_{n-1})-((\delta_{n}-\delta_{n-1})y_{n-1}+(\alpha_{n-1}-\alpha_{n})Ay_{n-1})\\
&&+\gamma\alpha_{n}(f(x_{n})-f(x_{n-1}))+
\gamma(\alpha_{n}-\alpha_{n-1})f(x_{n-1})\\
\end{eqnarray*}
\begin{eqnarray*}
&&+\delta_{n}x_{n}-\delta_{n}x_{n-1}+\delta_{n}x_{n-1}-\delta_{n-1}x_{n-1}\|\\
&\leq&(1-\delta_{n}-\alpha_{n}\bar{\gamma})\|y_{n}-y_{n-1}\|+|\delta_{n}-\delta_{n-1}|\|y_{n-1}\|\\
&&+|\alpha_{n}-\alpha_{n-1}|\|Ay_{n-1}\|+\gamma\alpha_{n}\alpha\|x_{n}-x_{n-1}\|\\
&&+\gamma|\alpha_{n}-\alpha_{n-1}|\|f(x_{n-1})\|+\delta_{n}\|x_{n}-x_{n-1}\|+|\delta_{n}-\delta_{n-1}|\|x_{n-1}\|\\
&\leq&(1-\delta_{n}-\alpha_{n}\bar{\gamma})[\|x_{n}-x_{n-1}\|+(1-\beta_{n})|\gamma_{n}-\gamma_{n-1}|\|x_{n-1}-W_{n-1}x_{n-1}\|\\
&&+(1-\beta_{n})(1-\gamma_{n})M_{1}\prod_{i=1}^{n-1}\gamma_{i}+(1-\beta_{n})M_{2}\prod_{i=1}^{n-1}\gamma_{i}\\
&&+|\beta_{n}-\beta_{n-1}|\|x_{n-1}-W_{n-1}z_{n-1}\|]+|\delta_{n}-\delta_{n-1}|\|y_{n-1}\|+|\alpha_{n}-\alpha_{n-1}|\|Ay_{n-1}\|\\
&&+\gamma\alpha_{n}\alpha\|x_{n}-x_{n-1}\|+\gamma|\alpha_{n}-\alpha_{n-1}|\|f(x_{n-1})\|+\delta_{n}\|x_{n}-x_{n-1}\|\\
&&+|\delta_{n}-\delta_{n-1}|\|x_{n-1}\|\\
&=&(1-\alpha_{n}\bar{\gamma})\|x_{n}-x_{n-1}\|+(1-\delta_{n}-\alpha_{n}\bar{\gamma})[(1-\beta_{n})|\gamma_{n}-\gamma_{n-1}|\|x_{n-1}-W_{n-1}x_{n-1}\|\\
&&+(1-\beta_{n})(1-\gamma_{n})M_{1}\prod_{i=1}^{n-1}\gamma_{i}+(1-\beta_{n})M_{2}\prod_{i=1}^{n-1}\gamma_{i}\\
&&+|\beta_{n}-\beta_{n-1}|\|x_{n-1}-W_{n-1}z_{n-1}\|]+|\alpha_{n}-\alpha_{n-1}|\|Ay_{n-1}\|+\gamma\alpha_{n}\alpha\|x_{n}-x_{n-1}\|\\
&&+\gamma|\alpha_{n}-\alpha_{n-1}|\|f(x_{n-1})\|+|\delta_{n}-\delta_{n-1}|\|y_{n-1}\|+|\delta_{n}-\delta_{n-1}|\|x_{n-1}\|\\
&\leq&(1-\alpha_{n}(\bar{\gamma}-\gamma\alpha))\|x_{n}-x_{n-1}\|\\
&&+(1-\delta_{n}-\alpha_{n}\bar{\gamma})[(1-\beta_{n})|\gamma_{n}-\gamma_{n-1}|\sup\{\|x_{n-1}+
\|W_{n-1}x_{n-1}\|\}\\
&&+(1-\beta_{n})\left((1-\gamma_{n})M_{1}\prod_{i=1}^{n-1}\gamma_{i}+M_{2}\prod_{i=1}^{n-1}\gamma_{i}\right)]\\
&&+|\alpha_{n}-\alpha_{n-1}|\sup\{\|Ay_{n-1}\|+\gamma f(x_{n-1})\}+|\delta_{n}-\delta_{n-1}|\sup\{\|y_{n-1}\|\\
&&+\|x_{n-1}\|\}+ |\beta_{n}-\beta_{n-1}|\sup\{\|x_{n-1}\|+\|W_{n-1}z_{n-1}\|\} .\\
\end{eqnarray*}
Now by lemma (\ref{1}) and $C_{3}, C_{4}, C_{5}$  we have $\|x_{n}-x_{n-1}\| \to 0 $.\\

\noindent On the other hand,
\begin{eqnarray*}
\|x_{n+1}-y_{n}\|&=&\|P_{C}[\alpha_{n}\gamma f(x_{n})+\delta_{n}x_{n}+((1-\delta_{n})I-\alpha_{n}A)y_{n}]-P_{C}(y_{n})\|\\
&\leq&\|\alpha_{n}\gamma
f(x_{n})+\delta_{n}x_{n}+((1-\delta_{n})I-\alpha_{n}A)y_{n}-y_{n}\|\\
&=&\|\alpha_{n}\gamma
f(x_{n})+\delta_{n}x_{n}-\delta_{n}x_{n+1}+\delta_{n}x_{n+1}+y_{n}-\delta_{n}y_{n}-y_{n}-\alpha_{n}Ay_{n}\|\\
&=&\|\alpha_{n}\gamma
f(x_{n})+\delta_{n}(x_{n}-x_{n+1})+\delta_{n}(x_{n+1}-y_{n})-\alpha_{n}Ay_{n}\|\\
&\leq&\alpha_{n}\|\gamma
f(x_{n})-Ay_{n}\|+\delta_{n}\|x_{n}-x_{n+1}\|+\delta_{n}\|x_{n+1}-y_{n}\|.\\
\end{eqnarray*}

\noindent So,
$\|x_{n+1}-y_{n}\|\leq\frac{\alpha_{n}}{(1-\delta_{n})}\|\gamma
f(x_{n})-Ay_{n}\|+\frac{\delta_{n}}{(1-\delta_{n})}\|x_{n}-x_{n+1}\|,$ which implies, $\|x_{n+1}-y_{n}\|\to 0.$\\

\noindent Also we have
$\|x_{n}-y_{n}\|\leq\|x_{n}-x_{n+1}\|+\|x_{n+1}-y_{n}\| ,$ which implies $\|x_{n}-y_{n}\|\to0$.\\

\noindent Notice that, \\
$\|z_{n}-x_{n}\|=\|\gamma_{n}x_{n}+(1-\gamma_{n}W_{n}x_{n}-x_{n}\|=\|(\gamma_{n}-1)x_{n}+(1-\gamma_{n})W_{n}x_{n}\|$ ,\\
and,\\
$\|y_{n}-W_{n}z_{n}\|=\|\beta_{n}x_{n}+(1-\beta_{n})W_{n}z_{n}-W_{n}z_{n}\|=\beta_{n}\|x_{n}-W_{n}z_{n}\|$ .\\

\noindent By two above equalities we have,
\begin{eqnarray*}
\|W_{n}x_{n}-x_{n}\|&\leq&\|x_{n}-y_{n}\|+\|y_{n}-W_{n}x_{n}\|\\
&\leq&\|x_{n}-y_{n}\|+\|y_{n}-W_{n}z_{n}\|+\|W_{n}z_{n}-W_{n}x_{n}\|\\
&\leq&\|x_{n}-y_{n}\|+\beta_{n}\|x_{n}-W_{n}x_{n}\|+\beta_{n}\|W_{n}x_{n}-W_{n}z_{n}\|+\|z_{n}-x_{n}\|\\
&\leq&\|x_{n}-y_{n}\|+\beta_{n}\|x_{n}-W_{n}x_{n}\|+(1+\beta_{n})\|z_{n}-x_{n}\|\\
&\leq&\|x_{n}-y_{n}\|+\beta_{n}\|x_{n}-W_{n}x_{n}\|+(1-\gamma_{n})(1+\beta_{n})\|W_{n}x_{n}-x_{n}\| .\\
\end{eqnarray*}
Therefore,\\
$[(1+\beta_{n})\gamma_{n}-2\beta_{n}]\|W_{n}x_{n}-x_{n}\|\leq\|x_{n}-y_{n}\|\to 0,$ so $\lim_{n\to\infty}\|W_{n}x_{n}-x_{n}\|=0$.\\

\noindent Furthermore we have,\\
$\|Wx_{n}-x_{n}\|\leq\|Wx_{n}-W_{n}x_{n}\|+\|W_{n}x_{n}-x_{n}\|,$ hence $\lim_{n\to\infty}\|Wx_{n}-x_{n}\|=0$.\\

\noindent We show that $\limsup_{n\to\infty}\langle\gamma
f(q)-Aq,x_{n}-q\rangle\leq0$, where $q=\lim_{t\to0}x_{t}$ and $x_{t}$ is the fixed point of the contraction $x\mapsto t\gamma f(x)+(I-tA)Wx$.\\

\noindent We have,\\
$\|x_{t}-x_{n}\|=\|(I-tA)(Wx_{t}-x_{n})+t(\gamma f(x_{t})-Ax_{n})\|$
and by lemma (\ref{r}),
\begin{eqnarray}
\|x_{t}-x_{n}\|^2&=&\|(I-tA)(Wx_{t}-x_{n})+t(\gamma f(x_{t})-Ax_{n})\|^2\nonumber\\
&\leq&(1-t\bar{\gamma})^2\|Wx_{t}-x_{n}\|^2+2t\langle\gamma f(x_{t})-Ax_{n},x_{t}-x_{n}\rangle\nonumber\\
&\leq&(1-2\bar{\gamma}t+(\bar{\gamma}t)^2)\|x_{t}-x_{n}\|^2+f_{n}(t)+2t\langle \gamma f(x_{t})-Ax_{t},x_{t}-x_{n}\rangle\nonumber\\
&&+2t\langle Ax_{t}-Ax_{n},x_{t}-x_{n}\rangle\ \nonumber\label{g},\\
\end{eqnarray}
where
$f_{n}(t)=\left(2\|x_{t}-x_{n}\|+\|x_{n}-Wx_{n}\|\right)\|x_{n}-Wx_{n}\|\to
0 \ , \ \ (as \ n\to \infty).$
Since $A$ is strongly positive linear mapping, so we have, \\
$\langle Ax_{t}-Ax_{n},x_{t}-x_{n}\rangle=\langle A(x_{t}-x_{n}),x_{t}-x_{n}\rangle\geq\bar{\gamma}\|x_{t}-x_{n}\|^2 $.\\

\noindent From(\ref{g}) we have,\\
\begin{eqnarray*}
2t\langle Ax_{t}-\gamma
f(x_{t}),x_{t}-x_{n}\rangle&\leq&(\bar{\gamma}^2
t^2-2\bar{\gamma}t)\|x_{t}-x_{n}\|^2+ f_{n}(t)+2t\langle Ax_{t}-Ax_{n},x_{t}-x_{n}\rangle\\
&\leq&\left(\bar{\gamma} t^2\right) \langle A(x_{t}-x_{n}),x_{t}-x_{n}\rangle+ f_{n}(t)+2t\langle A(x_{t}-x_{n}),x_{t}-x_{n}\rangle\\
&=&\bar{\gamma}t^2\langle A(x_{t}-x_{n}),x_{t}-x_{n}\rangle+f_{n}(t),\\
\end{eqnarray*}
which implies, $\langle Ax_{t}-\gamma f(x_{t}),x_{t}-x_{n}\rangle\leq\frac{\bar{\gamma}t}{2}\langle A(x_{t})-A(x_{n}),x_{t}-x_{n}\rangle+\frac{f_{n}(t)}{2t}$ .\\
Letting $n\to\infty$,\\
\begin{equation}\label{k}
\limsup\langle Ax_{t}-\gamma
f(x_{t}),x_{t}-x_{n}\rangle\leq\frac{t}{2}M_{3}   ,
\end{equation}
where $M_{3}$ is a constant such that, $\bar{\gamma}\langle
Ax_{t}-Ax_{n},x_{t}-x_{n}\rangle\leq M_{3}, \forall t\in (0,\min\{\|A\|^{-1},1\})$ and $n\geq 1$, taking $ t\to 0$, from (\ref{k}) we have, \\
\begin{equation}\label{l}
\limsup_{t\to0}\limsup_{n\to\infty}\langle Ax_{t}-\gamma
f(x_{t}),x_{t}-x_{n}\rangle\leq0.
\end{equation}
On the other hand we have,\\
\begin{eqnarray*}
\langle\gamma f(q)-Aq,x_{n}-q\rangle&=&\langle\gamma f(q)-Aq,x_{n}-q\rangle\\
&&-\langle\gamma f(q)-Aq,x_{n}-x_{t}\rangle+\langle\gamma f(q)-Aq,x_{n}-x_{t}\rangle\\
&&-\langle\gamma f(q)-Ax_{t},x_{n}-x_{t}\rangle+\langle\gamma f(q)-Ax_{t},x_{n}-x_{t}\rangle\\
&&-\langle\gamma f(x_{t})-Ax_{t},x_{n}-x_{t}\rangle+\langle\gamma f(x_{t})-Ax_{t},x_{n}-x_{t}\rangle.\\
\end{eqnarray*}
So,\\
 $\langle\gamma f(q)-Aq,x_{n}-q\rangle=\langle\gamma
f(q)-Aq,x_{t}-q\rangle+\langle Ax_{t}-Aq,x_{n}-x_{t}\rangle+ \langle\gamma f(q)-\gamma f(x_{t}),x_{n}-x_{t}\rangle+\langle\gamma f(x_{t})-Ax_{t},x_{n}-x_{t}\rangle$.\\

\noindent Hence,\\ $\limsup_{n\to\infty}\langle\gamma
f(q)-Aq,x_{n}-q\rangle\leq \|\gamma
f(q)-Aq\|\|x_{t}-q\|+\|A\|\|x_{t}-q\|\limsup_{n\to\infty}\|x_{n}-x_{t}\|$
$+\alpha\gamma\|q-x_{t}\|\limsup_{n\to\infty}\|x_{n}-x_{t}\|+\limsup_{n\to\infty}\langle\gamma f(x_{t})-Ax_{t},x_{n}-x_{t}\rangle$.\\

\noindent Therefore we have from (\ref{l}),\\
\begin{eqnarray*}
\limsup_{n\to\infty}\langle\gamma f(q)-Aq,x_{n}-q\rangle&=&\limsup_{t\to0}\limsup_{n\to\infty}\langle\gamma f(q)-Aq,x_{n}-q\rangle\\
&\leq&\limsup_{t\to0}\|\gamma f(q)-Aq\|\|x_{t}-q\|\\
&&+\limsup_{t\to0}\|A\|\|x_{t}-q\|\limsup_{n\to\infty}\|x_{n}-x_{t}\|\\
&&+\limsup_{t\to0}\gamma\alpha\|q-x_{t}\|\limsup_{n\to\infty}\|x_{n}-x_{t}\|\\
&&+\limsup_{t\to0}\limsup_{n\to\infty}\langle\gamma f(x_{t})-Ax_{t},x_{n}-x_{t}\rangle\leq0.\\
\end{eqnarray*}
Similarly,\\
\begin{eqnarray*}
\langle\gamma f(q)-Aq,y_{n}-q\rangle&=&\langle\gamma
f(q)-Aq,y_{n}-x_{n}\rangle+\langle\gamma f(q)-Aq,x_{n}-q\rangle\\
&\leq&\|\gamma f(q)-Aq\|\|y_{n}-x_{n}\|+\langle\gamma
f(q)-Aq,x_{n}-q\rangle,\\
\end{eqnarray*}
\noindent then, $\limsup_{n\to\infty}\langle\gamma
f(q)-Aq,y_{n}-q\rangle\leq0$.\\

\noindent Finally we prove that $x_{n}\to q$.\\
\begin{eqnarray*}
\|x_{n+1}-q\|^2&=&\|P_{C}[\alpha_{n}\gamma f(x_{n})+\delta_{n}x_{n}+((1-\delta_{n})I-\alpha_{n}A)y_{n}]-P_{C}(q)\|^2\\
&\leq&\|\alpha_{n}\gamma f(x_{n})+\delta_{n}x_{n}+((1-\delta_{n})I-\alpha_{n}A)y_{n}-q\|^2\\
&=&\|\alpha_{n}(\gamma f(x_{n})-Aq)+\delta_{n}(x_{n}-q)+((1-\delta_{n})I-\alpha_{n}A)(y_{n}-q)\|^2\\
&=&\|((1-\delta_{n})I-\alpha_{n}A)(y_{n}-q)+\delta_{n}(x_{n}-q)+\alpha_{n}(\gamma f(x_{n})-Aq)\|^2\\
&=&\|((1-\delta_{n})I-\alpha_{n}A)(y_{n}-q)+\delta_{n}(x_{n}-q)\|^2\\
&&+\alpha_{n}^2\|\gamma f(x_{n})-Aq\|^2+2\delta_{n}\alpha_{n}\langle x_{n}-q,\gamma f(x_{n})-Aq\rangle\\
&&+2\alpha_{n}\langle((1-\delta_{n})I-\alpha_{n}A)(y_{n}-q),\gamma f(x_{n})-Aq\rangle\\
\end{eqnarray*}
\begin{eqnarray*}
&\leq&[((1-\delta_{n})-\alpha_{n}\bar{\gamma})\|y_{n}-q\|+\delta_{n}\|x_{n}-q\|]^2\\
&&+\alpha_{n}^2\|\gamma f(x_{n})-Aq\|^2+2\delta_{n}\alpha_{n}\langle x_{n}-q,\gamma f(x_{n})-Aq\rangle\\
&&+2\alpha_{n}\langle((1-\delta_{n})I-\alpha_{n}A)(y_{n}-q),\gamma f(x_{n})-Aq\rangle\\
&=&[((1-\delta_{n})-\alpha_{n}\bar{\gamma})\|y_{n}-q\|+\delta_{n}\|x_{n}-q\|]^2\\
&&+\alpha_{n}^2\|\gamma f(x_{n})-Aq\|^2+2\delta_{n}\alpha_{n}\gamma\langle x_{n}-q,f(x_{n})-f(q)\rangle\\
&&+2\delta_{n}\alpha_{n}\langle x_{n}-q,\gamma f(q)-Aq\rangle+2(1-\delta_{n})\gamma\alpha_{n}\langle y_{n}-q,f(x_{n})-f(q)\rangle\\
&&+2(1-\delta_{n})\alpha_{n}\langle y_{n}-q,\gamma f(q)-Aq\rangle-2\alpha_{n}^2\langle A(y_{n}-q),\gamma f(q)-Aq\rangle\\
&\leq&[((1-\delta_{n})-\alpha_{n}\bar{\gamma})\|x_{n}-q\|+\delta_{n}\|x_{n}-q\|]^2\\
&&+\alpha_{n}^2\|\gamma f(x_{n})-Aq\|^2+2\delta_{n}\alpha_{n}\gamma\alpha\| x_{n}-q\|^2\\
&&+2\delta_{n}\alpha_{n}\langle x_{n}-q,\gamma f(q)-Aq\rangle+2(1-\delta_{n})\gamma\alpha_{n}\alpha\|x_{n}-q\|^2\\
&&+2(1-\delta_{n})\alpha_{n}\langle y_{n}-q,\gamma f(q)-Aq\rangle-2\alpha_{n}^2\langle A(y_{n}-q),\gamma f(q)-Aq\rangle\\
&=&[(1-\alpha_{n}\bar{\gamma})^2+2\delta_{n}\alpha_{n}\gamma\alpha+2(1-\delta_{n})\gamma\alpha_{n}\alpha]\|x_{n}-q\|^2\\
&&+\alpha_{n}^2\|\gamma f(x_{n})-Aq\|^2+2\delta_{n}\alpha_{n}\langle x_{n}-q,\gamma f(q)-Aq\rangle\\
&&+2(1-\delta_{n})\alpha_{n}\langle y_{n}-q,\gamma f(q)-Aq\rangle-2\alpha_{n}^2\langle A(y_{n}-q),\gamma f(q)-Aq\rangle\\
&\leq&[1-2(\bar{\gamma}-\alpha\gamma)\alpha_{n}]\|x_{n}-q\|^2+\bar{\gamma}^2\alpha_{n}^2\|x_{n}-q\|^2\\
&&+\alpha_{n}^2\|\gamma f(x_{n})-Aq\|^2+2\delta_{n}\alpha_{n}\langle x_{n}-q,\gamma f(q)-Aq\rangle\\
&&+2(1-\delta_{n})\alpha_{n}\langle y_{n}-q,\gamma f(q)-Aq\rangle+2\alpha_{n}^2\|A(y_{n}-q)\|\|\gamma f(q)-Aq\|\\
&=&[1-2(\bar{\gamma}-\alpha\gamma)\alpha_{n}]\|x_{n}-q\|^2+\alpha_{n}\{\alpha_{n}[\bar{\gamma}^2\|x_{n}-q\|^2 \\
&&+\|\gamma f(x_{n})-Aq\|^2+2\|A(y_{n}-q)\|\|\gamma f(q)-Aq\|]+2\delta_{n}\langle x_{n}-q,\gamma f(q)-Aq\rangle\\
&&+2(1-\delta_{n})\langle y_{n}-q,\gamma f(q)-Aq\rangle\}.\\
\end{eqnarray*}
Since $\{x_{n}\}$, $\{f(x_{n})\}$ and $\|y_{n}-p\|$ are bounded, we can take a constant $M_{4}>0$ such that,\\
$\bar{\gamma}^2\|x_{n}-q\|^2+\|\gamma f(x_{n})-Aq\|^2+2\|A(y_{n}-q)\|\|\gamma f(q)-Aq\|\leq M_{4} ,\ \forall \ n\geq \ 0$,\\
then it follows that, $\|x_{n+1}-q\|^2\leq
[1-2(\bar{\gamma}-\alpha\gamma)\alpha_{n}]\|x_{n}-q\|^2+\alpha_{n}\sigma_{n}$,
where,
$$\sigma_{n}=2\delta_{n}\langle x_{n}-q,\gamma f(q)-Aq\rangle+2(1-\delta_{n})\langle y_{n}-q,\gamma
f(q)-Aq\rangle+\alpha_{n}M_{4}.$$
Finally, we have $\limsup_{n\to\infty}\sigma_{n}\leq0$ and by lemma (\ref{1}) $x_{n}\to q$ .\ \ \ \ \ \ \ \ \ \ \ \ \ \  $\Box$\\

\par

\noindent Similar proof, shows that the followings composite iteration converges to $q\in F$, which solves variational inequality,\\
\begin{eqnarray}
x_{0}&=&x\in C \ chosen \ arbitrary \nonumber,\\
z_{n}&=&\lambda_{n}x_{n}+(1-\lambda_{n})K_{n}x_{n}\nonumber,\\
y_{n}&=&\beta_{n}x_{n}+(1-\beta_{n})K_{n}z_{n}\nonumber,\\
x_{n+1}&=&P_{C}[\alpha_{n}\gamma
f(x_{n})+\delta_{n}x_{n}+((1-\delta_{n})I-\alpha_{n}A)y_{n}]\nonumber\label{c}.\\
\end{eqnarray}

\begin{cor}\label{co}
Let $C$ be a closed convex subset of a real Hilbert space $H$. Let
$f$ be a contraction of $C$ into itself, let $A$ be a strongly
positive linear bounded operator with coefficient $\bar{\gamma}>0$
and $\{T_{i}:C \longrightarrow C\}$ be a finite family of
nonexpansive mappings. Assume that
$0<\gamma<\frac{\bar{\gamma}}{\alpha}$ and
$F=\bigcap_{i=1}^{N}F(T_{i})\neq\emptyset$. Let $x_{0} \in C$, given
that $\{\alpha_{n}\}$, $\{\beta_{n}\}$ and $\{\delta_{n}\}$ be
sequences in $[0,1]$ satisfying the following \
conditions:\\
$C_{1}: \lim_{n \rightarrow \infty}\alpha_{n}=0 \ ,
\sum_{n=1}^\infty\alpha_{n}=\infty$,\\
$C_{2} : 0<\liminf_{n \to\infty}\delta_{n}\leq\limsup_{n\to\infty}\delta_{n}<1$,\\
$C_{3}\colon \sum_{n=1}^{\infty}|\lambda_{n,i}-\lambda_{n-1,i}|<\infty$ , for all $i=1,2,\ldots,N$,\\
$C_{4}\colon \sum_{n=1}^{\infty}|\alpha_{n}-\alpha_{n-1}|<\infty$,\\
$C_{5}\colon \sum_{n=1}^{\infty}|\beta_{n}-\beta_{n-1}|<\infty$,\\
$C_{6}\colon (1+\beta_{n})\gamma_{n}-2\beta_{n}>d\quad $ for some $\quad  d\in (0,1)$.\\
If $\{x_{n}\}_{n=1}^{\infty}$ is the composite process defined by (\ref{c})
then the sequence $\{x_{n}\}_{n=1}^{\infty}$ converges strongly
to $q\in F$, which solves variational inequality $\langle\gamma f(q)-Aq,p-q\rangle\leq0$ ,$\forall p\in F$. \\

\end{cor}
\noindent If $\lambda_{n}=1$ and $\delta_{n}=0$ in corollary \ref{co}, then we get the result of Singthong and Suantai \cite{sing}.\\

\begin{cor}
Let $H$ be a Hilbert space, $C$ a closed convex subset of $H$. Let
$A$ be a strongly positive linear bounded operator with coefficient
$\bar{\gamma}\geq0$, and $f$ is a contraction. Let
$\{T_{i}\}_{i}^{N}$ be a finite family of nonexpansive mappings of
$C$ into itself and let $K_{n}$ be defined by (\ref{64}). Assume
that $0<\gamma<\frac{\bar{\gamma}}{\alpha}$ and $F=\bigcap_{i=1}^N
F(T_{i})\neq\emptyset$. Let $x_{1}\in C$, given that
$\{\alpha_{n}\}_{n=0}^\infty$ and $\{\beta_{n}\}_{n=0}^\infty$ are
sequences in $(0,1)$, and suppose that the following conditions are
satisfied:\\
$C_{1}\colon \alpha_{n}\rightarrow 0$, $\sum_{n=1}^\infty\alpha_{n}=\infty$,\\
$C_{2}\colon 0<\liminf_{n \rightarrow\infty}\beta_{n}\leq\limsup_{n\rightarrow\infty}\beta_{n}<1$,\\
$C_{3}\colon \sum_{n=1}^{\infty}|\gamma_{n,i}-\gamma_{n-1,i}|<\infty$ for all $i=1,2,\ldots,N$,\\
$C_{4}\colon \sum_{n=1}^{\infty}|\alpha_{n+1}-\alpha_{n}|<\infty$,\\
$C_{5}\colon \sum_{n=1}^{\infty}|\beta_{n+1}-\beta_{n}|<\infty$.\\
If $\{x_{n}\}_{n=1}^\infty$ is the composite process defined by
(\ref{va}) then the sequence $\{x_{n}\}$ converges strongly to $q\in
F$, which solves the variational inequality $\langle\gamma
f(q)-Aq,p-q\rangle\leq0$ ,$\forall p\in F$.\\
\end{cor}

\section{Simulation examples}In this section, we give three numerical examples to support the theoretical results. The iterations have been carried out on MATLAB 7.12. Here we recall
$r(n)=log_{10}\|x_{n+1}-x_{n}\|$ and
$\delta(n)=log_{10}\frac{\|x_{n}-x^{*}\|}{\|x^{*}\|}$ (i.e.
$\delta(n)$ is relative error), where $x^{*}$ is a fixed point of
$W_{n}$-mapping or $K$-mapping.\\
\noindent In the following, we assume $\gamma_{1}=\frac{1}{2}$,
$\gamma_{2}=\frac{1}{3}$, $\gamma_{3}=\frac{1}{4}$, and $x_{0}=3$.\\

\begin{table}[h]
   \center{\begin{tabular}{|c|c|c|c|c|}     \hline
        & $x^{*}$ & iteration & $T_{1}(x^{*})$ & $T_{2}(x^{*})$  \\
         \hline   $W_{n}$ mapping & $0.75290$ & $25$ & $0.6837577884$ & $0.7297090424$     \\
         \hline   $K$ mapping & $0.71491$ & $19$ & $0.6555494556$ & 0.7551522437 \\
\hline     \end{tabular} \caption{$T_{1}(x)=sin(x)$ and
$T_{2}(x)=cos(x)$.}}
 \end{table}

 \begin{table}[h]

  \center{    \begin{tabular}{|c|c|c|c|c|}     \hline
        & $x^{*}$ & iteration & $T_{1}(x^{*})$ & $T_{3}(x^{*})$  \\

         \hline   $W_{n}$ mapping & $0.0089628$ & $44834$ & $0.0089626800$ & $0.0089625600$     \\
         \hline
                  $K$ mapping & $ 0.0080118$ & $40066$ & $0.0080117142$ & $0.0080116285$ \\
\hline

                              \end{tabular}
                                \caption{$T_{1}(x)=sin(x)$ and $T_{3}(x)=tan^{-1}(x)$.}
}
                                   \end{table}

\begin{table}[h]

      \center{\begin{tabular}{|c|c|c|c|c|c|}     \hline
        & $x^{*}$ & iteration & $T_{1}(x^{*})$ & $T_{2}(x^{*})$ & $T_{3}(x^{*})$  \\

         \hline   $W_{n}$ mapping & $0.59403$ & $85$ & $0.5597051868$ & $0.8286918026 $ & $0.5360182305$  \\
         \hline
                  $K$ mapping & $0.67735$ & $18$ & $0.6267302508$ & $0.7792362880$ & $0.5953623347$ \\
\hline

                              \end{tabular}
                                \caption{$T_{1}(x)=sin(x)$, $T_{2}(x)=cos(x)$ and $T_{3}(x)=tan^{-1}(x)$.}
}
                                   \end{table}
\begin{figure}\label{fig1}
\center\includegraphics[width=6cm]{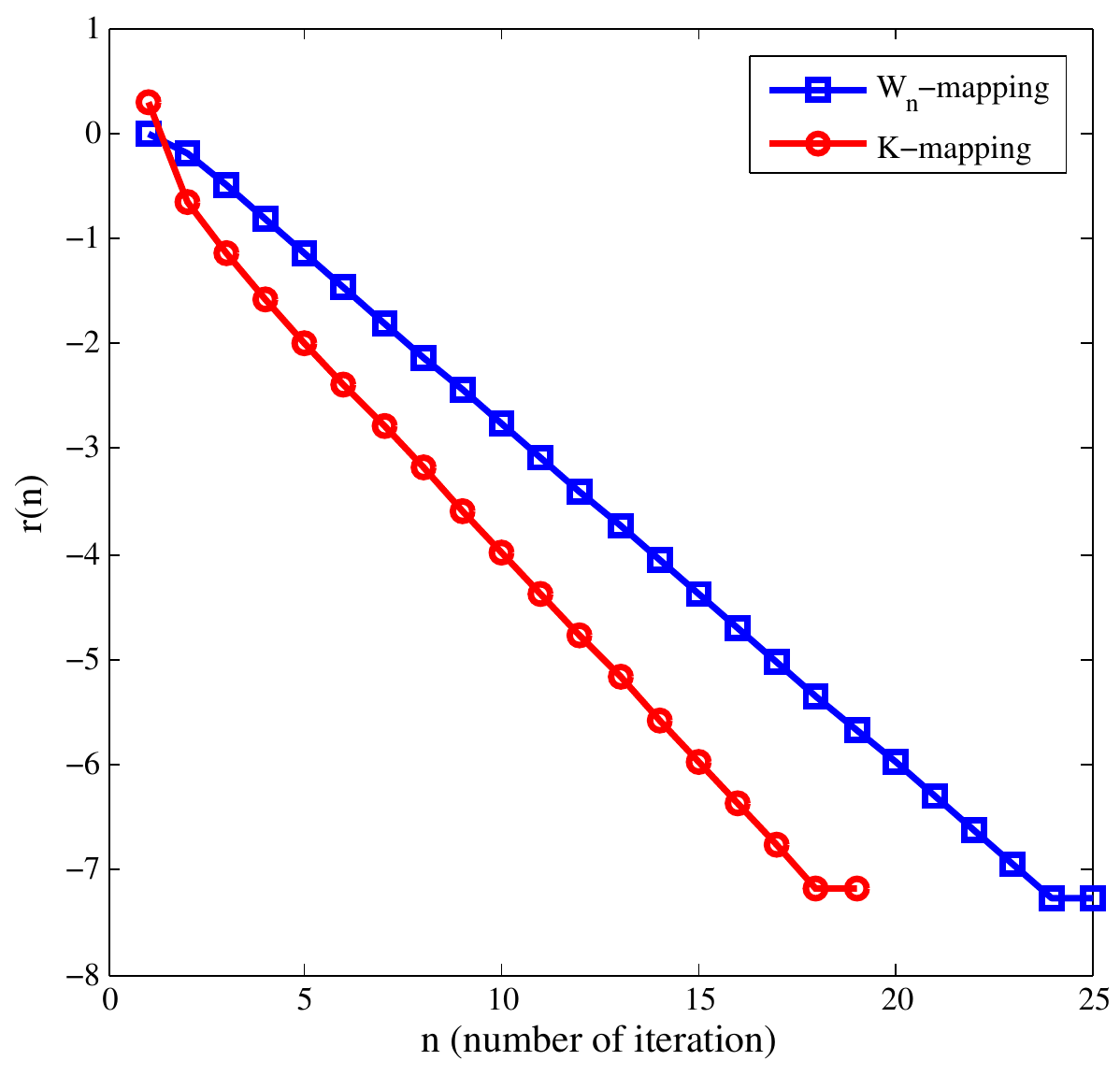}
\includegraphics[width=6.1cm]{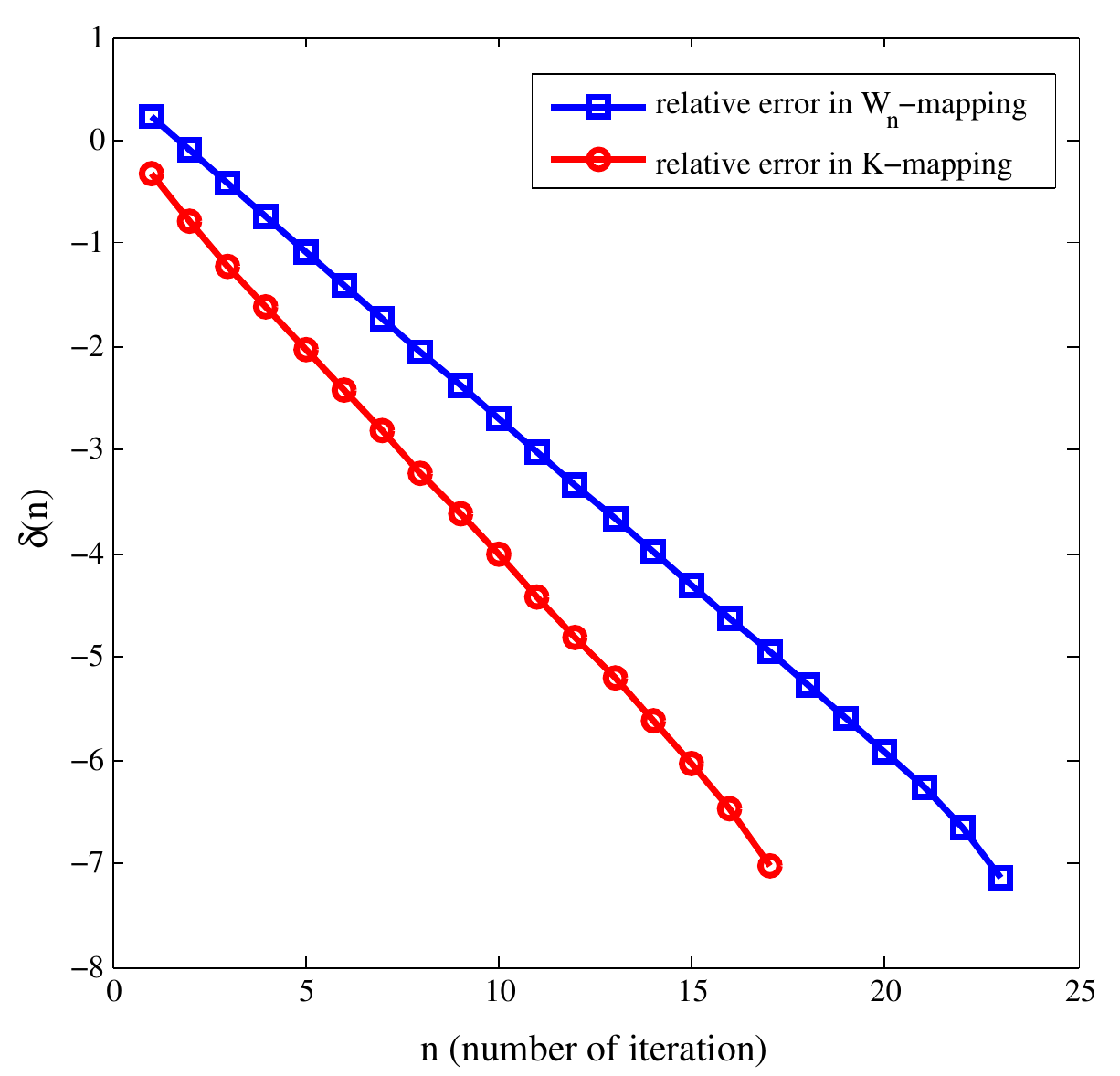}
\caption{\small{The results obtained for $T_{1}$ and $T_{2}$.}}
\end{figure}
\begin{figure}\label{fig2}
\center\includegraphics[width=6cm]{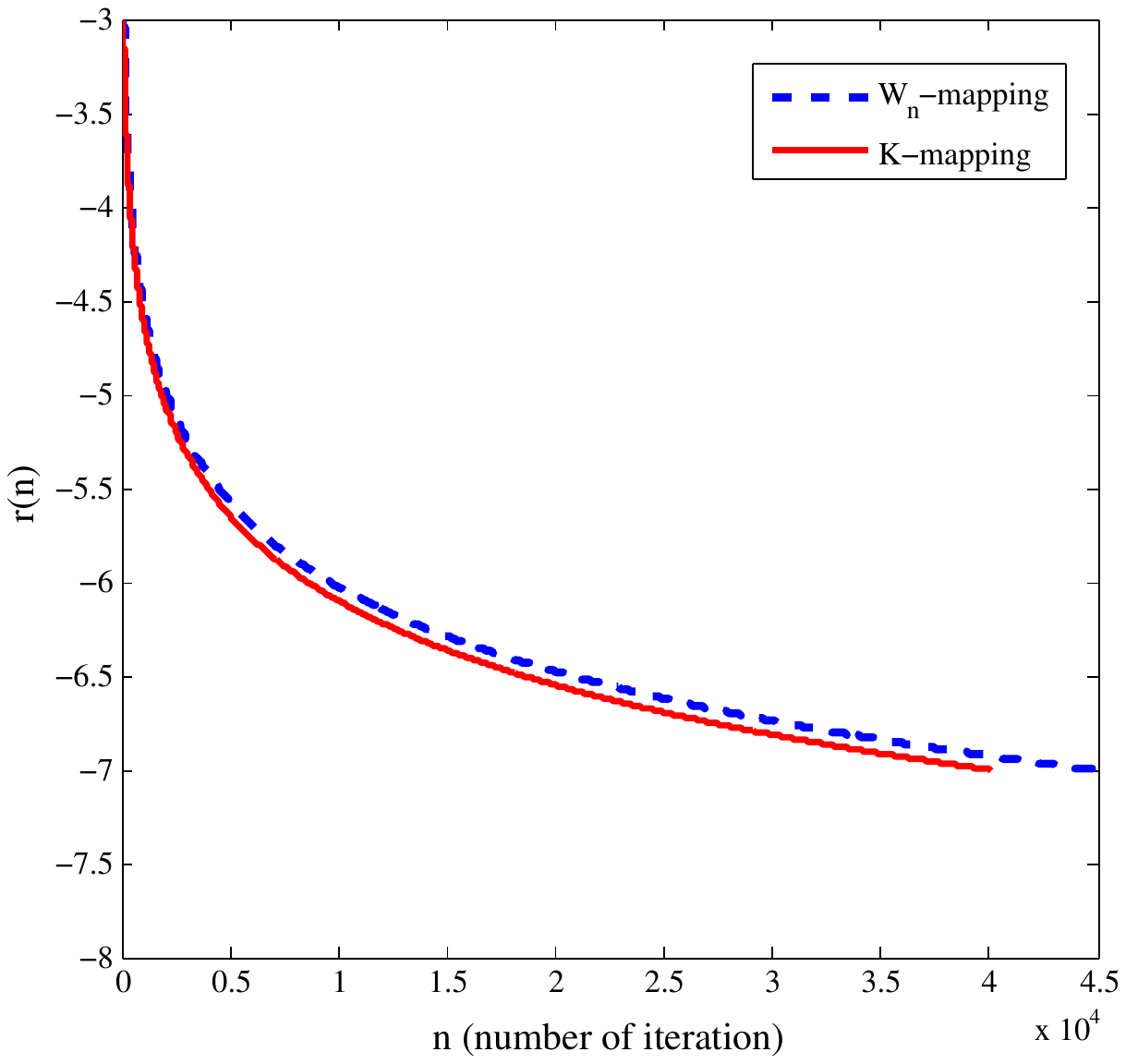}
\includegraphics[width=5.9cm]{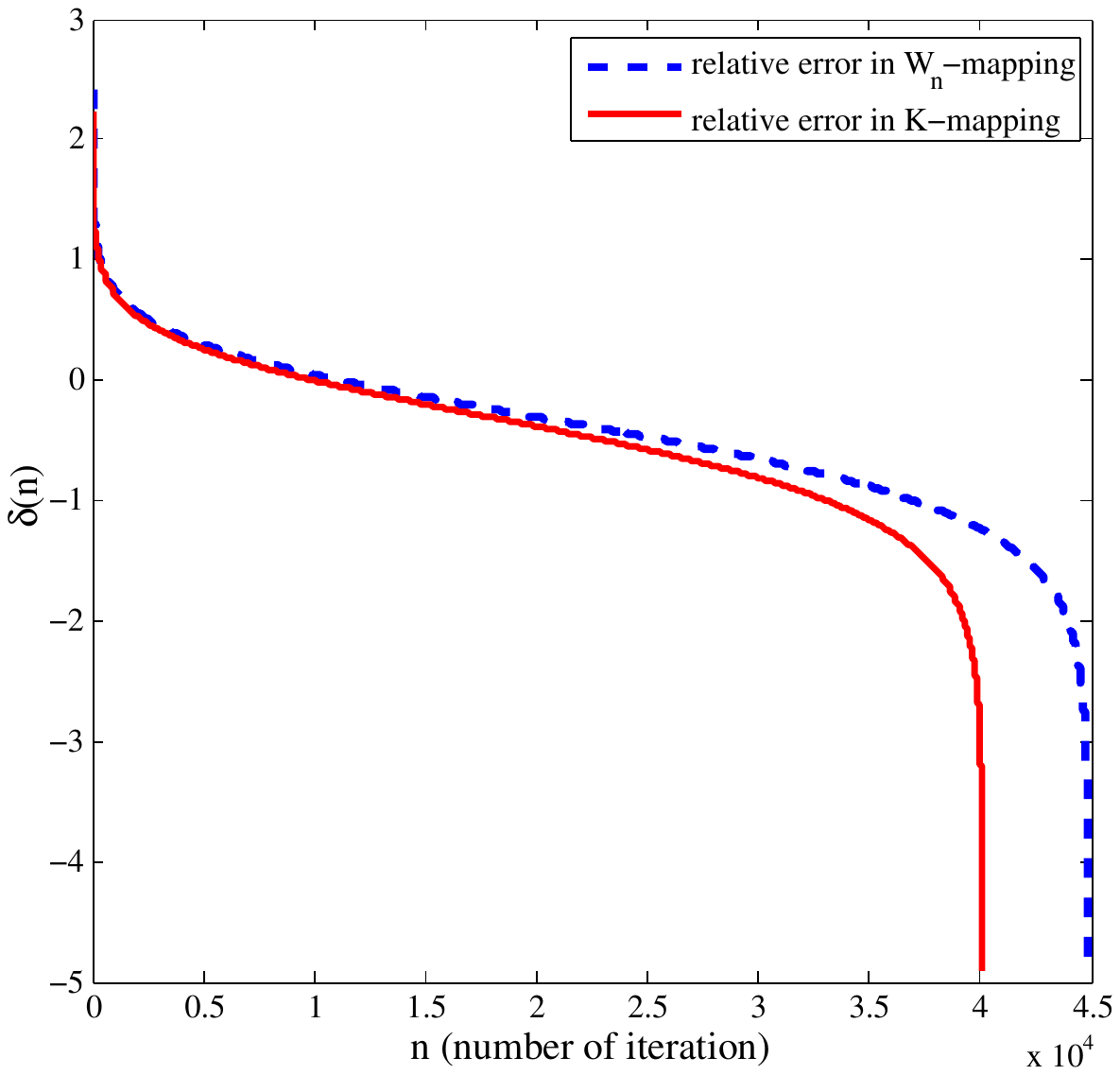}
\caption{\small{The results obtained for $T_{1}$ and $T_{3}$.}}
\end{figure}
\begin{figure}\label{fig3}
\center\includegraphics[width=6cm]{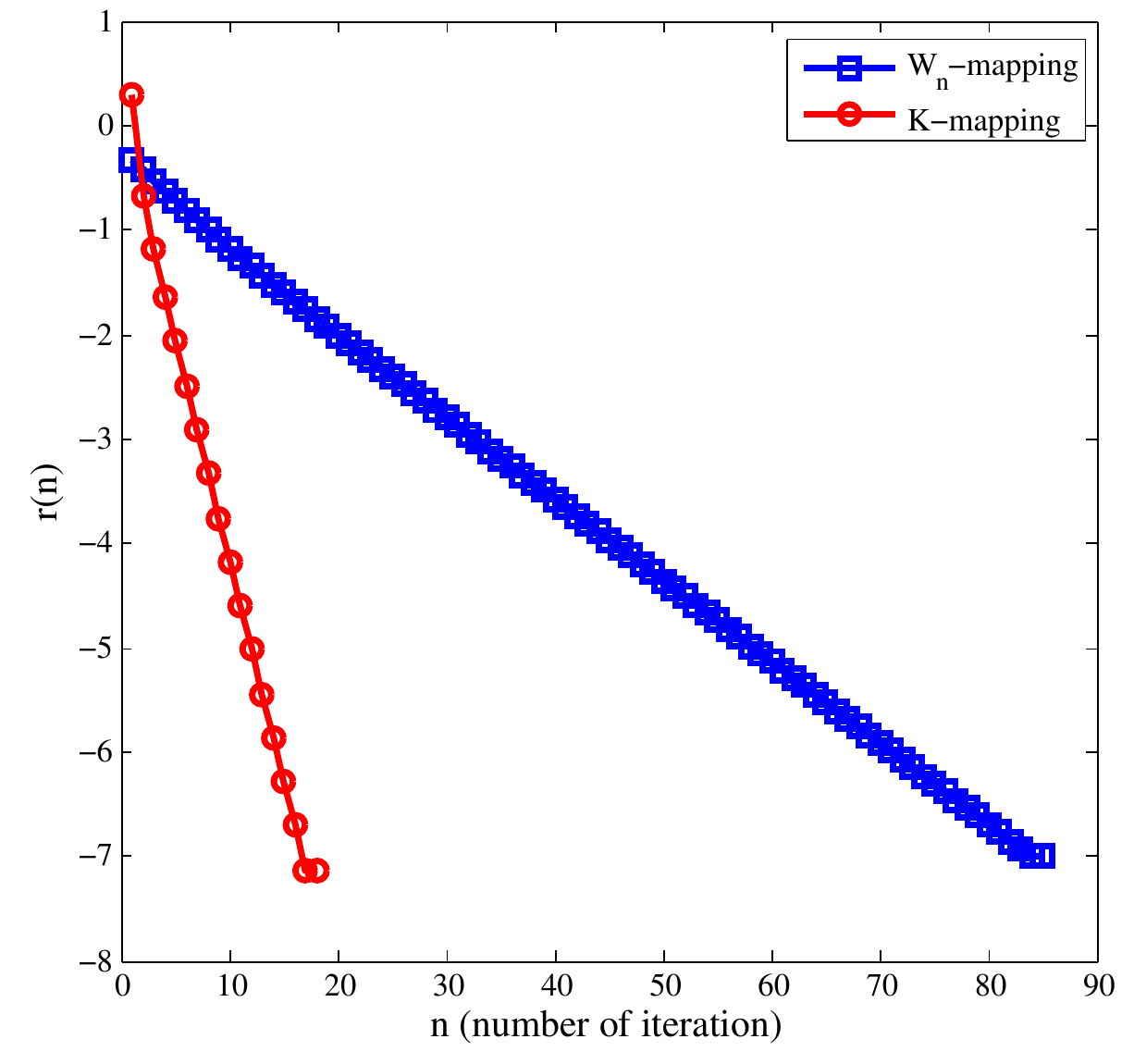}
\includegraphics[width=6cm]{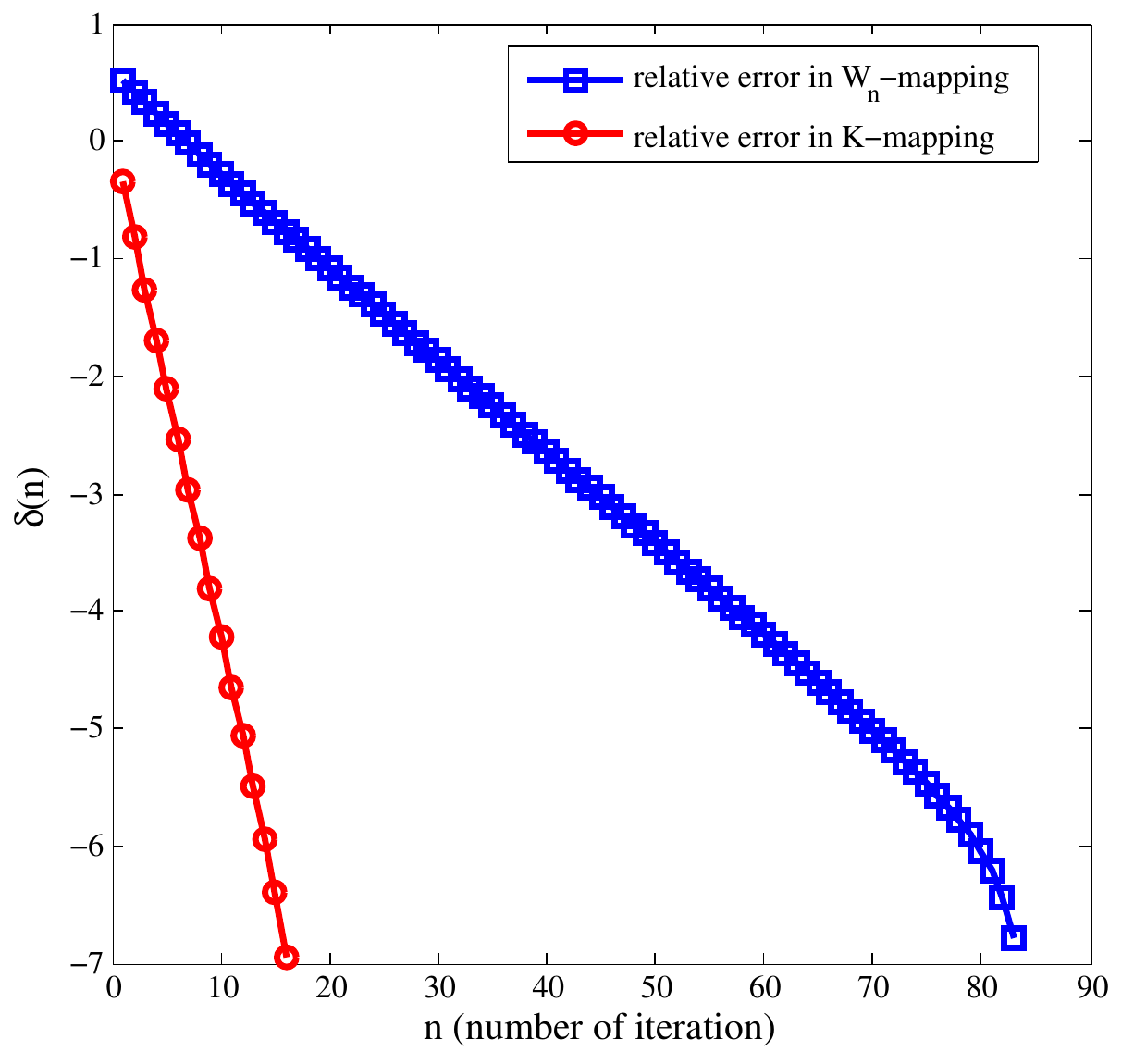}
\caption{\small{The results obtained for $T_{1}$, $T_{2}$ and
$T_{3}$.}}
\end{figure}

\section{Conclusion} Finding the fixed point of nonexpansive mappings and variational inequalities is so important in many fields. In this paper, we have
constructed an iterative algorithm for finding a common fixed point
of an infinite family of nonexpansive mappings and a solution of
certain variational inequality. Finally, some numerical examples
were presented to support the theoretical results of this paper.
Moreover, these examples compare the error and speed of convergence
of $W_{n}$-mapping and $K$-mapping.

\end{document}